\documentclass[a4paper]{article}
\usepackage[latin1]{inputenc}
\usepackage{fullpage}
\usepackage{amsmath,amsthm,amsfonts}
\usepackage{xspace,url}
\usepackage[pdftex,colorlinks=true]{hyperref}
\usepackage[numbers]{natbib}
\renewcommand{\v}[1]{\boldsymbol{#1}} % Vectors and bivectors.
\renewcommand{\i}{\ensuremath{\v{i}}\xspace}
\renewcommand{\j}{\ensuremath{\v{j}}\xspace}
\renewcommand{\k}{\ensuremath{\v{k}}\xspace}
\let\oldmu=\mu
\renewcommand{\mu}{\ensuremath{\v{\oldmu}}\xspace}
  \newcommand{\R}{\ensuremath{\mathbb{R}}\xspace}
  \newcommand{\C}{\ensuremath{\mathbb{C}}\xspace}
\newtheorem{lemma}{Lemma}
\newtheorem{theorem}{Theorem}
\title{Quaternion polar representation with a complex modulus and complex argument inspired by the Cayley-Dickson form}
\author{Stephen~J.~Sangwine\thanks{Department of Computing and Electronic Systems,
                 University of Essex, Wivenhoe Park, Colchester, CO4 3SQ,
                 United Kingdom.
                 Email:~\texttt{S.Sangwine@IEEE.org}}
       \and Nicolas~Le~Bihan\thanks{Département Images et Signal,
                 Gipsa-lab UMR 5216,
                 961 Rue de la Houille Blanche,
                 Domaine Universitaire, BP 46,
                 38402 Saint Martin d'Hères,
                 Cedex, France.
                 Email:~\texttt{nicolas.le-bihan@gipsa-lab.inpg.fr}}}
\begin{document}
\maketitle
\begin{abstract}
We present a new polar representation of quaternions inspired by the Cayley-Dickson
representation. In this new polar representation, a quaternion is represented by a
pair of complex numbers as in the Cayley-Dickson form, but here these two complex
numbers are a complex `modulus' and a complex `argument'. As in the Cayley-Dickson
form, the two complex numbers are in the same complex plane (using the same complex
root of $-1$), but the complex phase is multiplied by a different complex root of $-1$
in the exponential function. We show how to calculate
the amplitude and phase from an arbitrary quaternion in Cartesian form.
\end{abstract}
\section{Introduction}
\label{intro}
It is well-known that the complex exponential $e^{\i\theta} = \cos\theta+\i\sin\theta$
generalises to quaternions by replacing \i by any unit pure quaternion \mu, since any
unit pure quaternion is a root of $-1$
\cite[pp.\,203, 209]{Hamiltonpapers:V3:7}\nocite{Hamilton:1848}\cite[\S167, p179]{Hamilton:1853}.
Hence, any quaternion may be represented in the
polar form $q = |q|e^{\mu\theta}$ where $\theta$ is a real angle \cite[§2.3]{Ward:1997}\cite[§12.7]{Altmann:1986}.
A small difference between the complex and quaternion cases is that in the quaternion case
the argument $\theta$ is conventionally confined to the interval $[0,\pi]$. This is
because the modulus of the vector part of the quaternion is always taken to be positive
(there is no convenient way to define an orientation in 3-space which would permit the sign
of the vector part to be determined).

The classic polar form just discussed is not the only possibility for representing a
quaternion. There are representations based on Euler angles using three exponentials in
which the arguments of the exponentials are directly linked to the Euler angles.
Bülow, for example, in  \cite[Theorem 2.3]{Bulow:1999} and \cite[§V.A]{BulowSommer:2001}
quotes the following formula, and gives an algorithm for finding the three angles:
\[
q = |q|e^{\i\phi}e^{\k\psi}e^{\j\theta}
\]
See also \cite[§12.11]{Altmann:1986} for a discussion of Euler angles and quaternions.

The polar form we present in this paper is different to both of the preceding forms,
in that it represents a quaternion using a single exponential with a complex argument
and a complex `modulus'.

\section{Polar form}
We now show that every quaternion admits the following polar form, in addition to those
presented in the previous section:
\[
q = Ae^{B\j}
\]
where $A = a + b\i$ and $B = c + d\i$ are complex and $a$, $b$, $c$ and $d$ are real.
Before proceeding, we note the analogy between
this polar form, based on two complex numbers, and the Cayley-Dickson form of a quaternion
$q = (w + x\i) + (y + z\i)\j = w + x\i + y\j + z\k$ which is also based on two complex numbers. In the Cayley-Dickson
form the first complex number has real and imaginary parts which are simply the first two parts
of the quaternion in Cartesian form, and the second complex number has real and imaginary parts
which are simply the third and fourth parts of the quaternion in Cartesian form. There is a
subtlety to the Cayley-Dickson form which we must mention before proceeding, otherwise what
follows may not be clear. The subtlety lies in the difference between a complex number such as
$a + b\i$ and a quaternion $a + b\i + 0\j + 0\k$. In this paper, as in the classic Cayley-Dickson
form, the complex number is treated as a degenerate quaternion and we handle complex numbers
algebraically as quaternions in which two components happen to be zero. The usual rules of
quaternion algebra therefore apply, including the rules $\i^2 = \j^2 = \k^2 = \i\j\k = -1$. 

The polar form presented in this paper is very subtle. It shares with the Cayley-Dickson form
the idea of a construction based on two complex numbers in the same complex plane with \i as the
root of $-1$, the second of these (the argument of the exponential) being multiplied on the right
by \j, just as in the Cayley-Dickson form. Unlike the Cayley-Dickson form, the parameters of the
polar form presented here are not trivially obtainable from either the Cartesian form of a quaternion
or the classic polar form. 

In what follows, we first need to prove that $e^{B\j}$ is a special type of quaternion,
with only three components: $\alpha + \gamma\j + \delta\k$, that is, it has a zero coefficient of $\i$.

\begin{lemma}
\label{lemma1}
Given an arbitrary quaternion in the form $p = (c + d\i)\j = c\j + d\k$, its exponential
is given by:
\[
e^p = \cos|p| + \frac{p}{|p|}\sin|p| = \alpha + \gamma\j + \delta\k
\]
\end{lemma}
\begin{proof}
Represent $p$ by the product of its modulus and a unit quaternion obtained by dividing by the modulus:
\[
p = |p|\frac{p}{|p|}
\]
Comparison with the classic polar form of a quaternion given in §\ref{intro}, namely
$e^{\mu\theta} = \cos\theta + \mu\sin\theta$, where $\mu$
is a unit pure quaternion, is sufficient to prove the result, since $p/|p|$
is a unit pure quaternion which may be identified with $\mu$, and $|p|$ is real and may be
identified with $\theta$. Thus $\alpha = \cos|p|$, $\gamma = (c/|p|)\sin|p|$ and $\delta = (d/|p|)\sin|p|$.
\end{proof}
Now we present a polar form as previously outlined with a complex `modulus' and a complex
argument, and we show how to express an arbitrary quaternion in this form, that is, how to
find the complex `modulus' and argument.
\begin{theorem}
\label{thetheorem}
Every quaternion $q = w + x\i + y\j + z\k$, where $w, x, y, z\in\R$ can be expressed in the form
$q = Ae^{B\j}$, where $A, B\in\C$, specifically $A = a + b\i$ and $B = c + d\i$.
\end{theorem}
\begin{proof}
The proof is by demonstrating that we can find $A$ and $B$ without any constraints on $q$. Clearly,
if $q=0$, $A$ must be zero, and the value of $B$ is unimportant.
Otherwise, we note that if we put the
complex `modulus' $A$ in polar form as $A = |A|e^{\theta\i}$, we can immediately see that $|A|=|q|$.
Therefore in what follows we may assume, without loss of generality, that we are working with
quaternions $q$ of unit modulus, and that therefore $A = a + b\i$ is of unit modulus.
Therefore we have, using the result in Lemma \ref{lemma1}:
\[
q = Ae^{B\j} = (a + b\i)(\alpha + \gamma\j + \delta\k) = a\alpha + b\alpha\i + (a\gamma - b\delta)\j + (a\delta + b\gamma)\k
\]
Now it is possible to see that we can obtain $a + b\i$ quite trivially, provided $\alpha\ne 0$.
Assuming $\alpha\ne 0$ for the moment, and since $a + b\i$ has unit
modulus, all we need is to take a complex number formed from $a\alpha$ and $b\alpha$ (the first and second
components of $q$, which are the scalar part and the coefficient of \i), and normalise this complex number
by dividing by its modulus.
However, there is a small problem: there is an ambiguity of sign. We do not know the sign of
$\alpha$, and this means we have not determined $a + b\i$ completely since $-(a + b\i)$ is equally valid. It
is clear that this ambiguity is insurmountable, since negating $a + b\i$ can be compensated for by negating
$\alpha + \gamma\j + \delta\k$.

In the case where $\alpha=0$, but $|q|>0$, the first and second components of $q$ will be zero.
Thus, $q$ is of the form $\gamma\j+\delta\k$ and we may take $a=1$ and $b=0$ (thus $A=1$).
The factorization here is not unique, and in §\ref{numeric} we present an example to show
this (we could have chosen $a=0$ and $b=1$, and thus $A=\i$).

Once we have obtained $a + b\i$ it is simple to find $\alpha + \gamma\j + \delta\k$ by dividing $q$ on the
left by $a + b\i$, or equivalently, by multiplying on the left by the inverse\footnote{The inverse of
$a+b\i$ is $a - b\i$ if $a + b\i$ has unit modulus (as assumed) since the inverse is the conjugate
divided by the modulus.} of $a + b\i$.

To obtain $B\j$ we can use the quaternion $\log$ function:
\[
B\j = \log(\alpha + \gamma\j + \delta\k)
\]
Numerically, this is a simple process, since the quaternion logarithm can be computed using a complex
logarithm operating on a complex number isomorphic to the quaternion. The resulting complex value can
then be used trivially to construct a quaternion with vector part oriented in the same direction as
the original quaternion\footnote{See, for example, the file \texttt{log.m} in our Matlab® quaternion
toolbox \cite{QTFM}}.

It should be clear from Lemma \ref{lemma1} that $B\j$ must be of the form $c\j+d\k$ and therefore the
real and imaginary parts of $B$ are trivially obtained from the second and third imaginary components
of $B\j$ (the coefficients of \j and \k). 
\end{proof}
\subsection{Sign ambiguity}
We noted above that the sign of $A$ has an ambiguity which can be compensated by an alteration to $B$
but we did not discuss what form this alteration would take. We merely pointed out that negating $A$
required $\alpha + \gamma\j + \delta\k$ to be negated.

It is easy to see that negation of $\alpha + \gamma\j + \delta\k$ may be accomplished by \emph{adding} $\pi$ to the
modulus of $B$. From Lemma \ref{lemma1} adding $\pi$ to $|B|$ means that the angles in the
cosine and sine are augmented by $\pi$ which negates their results. It is not so easy to see that
the converse is not the only possibility. In fact alternative solutions exist that modify $B$ to
negate $\alpha + \gamma\j + \delta\k$, and the algorithm presented above for finding $B$ using a quaternion
logarithm will yield these solutions. For the present paper we confine our discussion of this
issue to some numerical examples.

\section{Numerical examples}
\label{numeric}
Let $q = \frac{1}{2}(1 + \i + \j + \k)$. Then $|q| = 1$ and
$a\alpha = b\alpha = \frac{1}{2}$ which gives $A = \frac{1}{\sqrt{2}} + \frac{1}{\sqrt{2}}\i$.
Multiplying $q$ on the left by the conjugate of $A$ gives $\frac{1}{\sqrt{2}} + \frac{1}{\sqrt{2}}\j$
and the quaternion logarithm gives $B\j = \pi\j/4$ (therefore $B$ is real in this instance).
To verify the result:
\[
(\frac{1}{\sqrt{2}} + \frac{1}{\sqrt{2}}\i)\exp\frac{\pi}{4}\j =
(\frac{1}{\sqrt{2}} + \frac{1}{\sqrt{2}}\i)(\cos\frac{\pi}{4} + \j\sin\frac{\pi}{4}) =
(\frac{1}{\sqrt{2}} + \frac{1}{\sqrt{2}}\i)(\frac{1}{\sqrt{2}} + \j\frac{1}{\sqrt{2}}) =
\frac{1}{2}(1 + \i + \j + \k)
\]

For a second example, we take $q = 1 + 2\i + 3\j + 4\k$. In this case we state the result
obtained numerically which is $A = 2.4495 + 4.8990\i$ and $B = 1.1317 - 0.2058\i$. The
modulus of $B = 1.1503$. Hence (apart from rounding errors caused by the limited precision
of the intermediate results shown here):
\begin{equation*}
\begin{split}
(2.4495 + 4.8990\i)\exp(1.1317\j - 0.2058\k) =\\
(2.4495 + 4.8990\i)(\cos 1.1503 + (0.9839\j - 0.1789\k)\sin 1.1503) =\\
(2.4495 + 4.8990\i)(0.40825 + 0.89815\j - 0.1633\k) = 1 + 2\i + 3\j + 4\k
\end{split}
\end{equation*}

Finally, we show an example where the value $\alpha$ in Theorem \ref{thetheorem} is
zero, and the factorization is not unique. Consider the factorization of the quaternion
with value $\k$. The following two factorizations will yield this value. Firstly $A=\i$
and $B=\frac{\pi}{2}$:
\[
\i\exp\left(\frac{\pi}{2}\j\right) = \i\left(\cos\frac{\pi}{2} + \j\sin\frac{\pi}{2}\right) = \i(0 + \j) = \k
\]
Secondly, $A=1$ and $B=\frac{\pi}{2}\i$:
\[
1\exp\left(\frac{\pi}{2}\i\j\right) = \cos\frac{\pi}{2} + \i\j\sin\frac{\pi}{2} = 0 + \i\j = \k
\]

%\bibliographystyle{plainnat}
%\bibliography{IEEEabrv,sangwine,quaternions,maths,clifford,online,image_processing}

\end{document}